\documentclass[10pt]{article}

\usepackage[utf8]{inputenc}
\usepackage{amscd,amssymb}
\usepackage{amsthm,amsmath}
\usepackage{graphicx}
\usepackage{enumerate,amssymb,setspace}
\usepackage{amsfonts,amsthm,amscd}
\usepackage{mathtools}
\usepackage{mathtools}
\usepackage[all]{xy}
\usepackage{hyperref}

\usepackage[margin=2cm]{geometry}
\usepackage{caption}
\usepackage[]{subcaption}

\DeclareCaptionLabelFormat{opening}{}

\newtheorem{theorem}{Theorem}
\newtheorem{lemma}[theorem]{Lemma}
\newtheorem{prop}[theorem]{Proposition}
\newtheorem{corollary}[theorem]{Corollary}
\newtheorem{conjecture}[theorem]{{Conjecture}}
\newtheorem{example}[theorem]{{Example}}
\newtheorem{problem}[theorem]{{Problem}}

\newtheorem{definition}[theorem]{{Definition}}
\newtheorem{claim}[theorem]{{Claim}}

\def\bclaim{\begin{claim}}
	\def\eclaim{\end{claim}}
\def\bdefin{\begin{definition}}
	\def\edefin{\end{definition}}
\def\bcor{\begin{corollary}}
	\def\ecor{\end{corollary}}
\def\bthm{\begin{theorem}}
	\def\ethm{\end{theorem}}
\def\bconj{\begin{conjecture}}
	\def\econj{\end{conjecture}}
\def\blem{\begin{lemma}}
	\def\elem{\end{lemma}}
\def\blemma{\begin{lemma}}
	\def\elemma{\end{lemma}}
\def\bprop{\begin{prop}}
	\def\eprop{\end{prop}}
\def\bremark{\begin{remark}}
	\def\eremark{\end{remark}}
\def\bprob{\begin{problem}}
	\def\eprob{\end{problem}}

\theoremstyle{remark}
\newtheorem{remark}[theorem]{Remark}

\newcommand{\Pic}{\mathrm{Pic}}

\newcommand{\rk}{\mathrm{rk}}

\makeatletter\@addtoreset{equation}{section} \makeatother

\def\be{\beta}

\def\h#1{\hbox{#1}}

\def\Pic{\operatorname{Pic}}

\newcommand{\QQ}{\mathbb{Q}} 
 \newcommand{\NN}{{\mathbb N}}

\def\beq{\begin{equation}}
	\def\eeq{\end{equation}}
\def\bpf{\begin{proof}}
	\def\epf{\end{proof}}
\def\bremark{\begin{remark}}
	\def\eremark{\end{remark}}

\def\eaeq{\end{aligned}}
\def\baeq{\begin{aligned}}

\newcommand\blfootnote[1]{%
	\begingroup
	\renewcommand\thefootnote{}\footnote{#1}%
	\addtocounter{footnote}{-1}%
	\endgroup
}

\title{Two remarks on asymptotically log Fano pairs\texttt{\blfootnote{MSC subject codes: 14J45, 14J26 (primary), 14J10, 14E05 (secondary).
			Keywords: asymptotically log Fano varieties, asymptotically log Del Pezzo surfaces, body of ample angles. 
		}
}	}

\author{Jesus Martinez-Garcia}
\begin{document}

	\maketitle
%
	

	Asymptotically log Fano pairs were introduced in \cite{ALF} by Cheltsov and Rubinstein, generalising a definition of \cite{Maeda}. They have received attention in the last decade within the theory of K-stability, e.g. see \cite{Fuj3, MGPZ}, as they approximate log Calabi Yau pairs while staying in the log Fano setting. In this note, written for the ZAG Proceedings, we summarise our talk on 1st September 2020 (Day of Knowledge) given on Zoom during the 24-hour ZAG Marathon. In the talk, we reported on joint work with P. Cascini and Y. Rubinstein \cite{ALF-ample-angles}, on the classification of the two-dimensional case, known as asymptotically log del Pezzo pairs.
	
	Let $X$ be a a smooth projective complex variety $X$ of dimension $n$ and $r\geqslant 1$. Let $D_i$, $i=1,\ldots,r$ be distinct irreducible hypersurfaces on $X$. The pair $(X,D=\sum_{i=1}^rD_i)$ is an \emph{asymptotically log Fano (ALF) variety} if there exists a sequence $\be(j)=(\be_1(j),\ldots,\be_r(j)) \in(0,1]^r\cap\QQ^r, \,j\in\NN$ converging to the origin such that
	$$-K_{(X,D)}^{\beta(j)}\coloneqq -K_X-\sum_{i=1}^r(1-\be_i(j))D_i$$
	is ample for all $j$. The simplest example is to take a Fano variety $X$ (e.g. $X=\mathbb P^n$) and a smooth $D\in |-K_X|$. The moduli space of the latter pairs in dimension $2$ is considered in \cite{MGPZ}. It is natural to ask the specific set of coefficients for which the ample condition holds. The \emph{body of ample angles} of $(X,D)$ is the set
	$$
				\h{\rm AA}(X,D):=\Big\{\be=(\be_1,\ldots,\be_r)\in(0,1)^r\,:\, 
				\h{$-K_X-\sum_{i=1}^r(1-\be_i)D_i$ is ample} \Big\},
	$$
	which is naturally convex since ampleness is a convex condition. The pair $(X,D)$ is \emph{strongly asymptotically log Fano} if there is some open ball $B_\varepsilon$, $\varepsilon>0$ centred at the origin such that $B_\varepsilon\cap (0,1]^r\subseteq \rm{AA}(X,D)$. It is worth noting this is not the definition given in \cite{ALF} but it is equivalent to it.

	All strongly ALF pairs are ALF, but the converse is not necessarily true. Indeed, when $r=1$ or $n=1$, the two statements are equivalent, but they differ when $r>2$ and $n\geqslant 2$, as the following example illustrates.
	
	\begin{example}
		\label{exa:ALDP.1.n}
	Consider the $n$-th Hirzebruch surface $X=\mathbb F_n=\mathbb P_{\mathbb {P}^1}(\mathcal O_{\mathbb P^1}\oplus \mathcal O_{\mathbb P^1}(n))$, $n>0$, let $D_1=Z_n$ be the unique curve in $X$ with negative self-intersection and $F$ be the class of a fibre of the fibration $\mathbb F_n\rightarrow \mathbb P^1$. Let $D_2\in |Z_n+(n+2)F|$. A divisor numerically equivalent to $aZ_n+bF$ is ample if and only if $b>na>0$. Thus,  $-K_{(X,D)}^{\beta(j)}$ is ample if and only if $2\beta_2>n\beta_1$, making $(X,D)$ an ALF, non-strongly ALF pair. Indeed, a quick computation shows that 
	$${\rm AA}(S,C)=\{(\be_1,\be_2)\in(0,1]^2\,:\,-n\be_1+2\be_2>0\}.$$
	
	\end{example}
	
%
	The first thing one notices from the example above is that the body of ample angles is, in this example, a polyhedron. Many invariants related to positivity of Fano varieties (e.g. the Mori cone) are polyhedral, so one wonders if this will be the case for the body of ample angles and it is indeed the case. In fact, more is true:
	\begin{theorem}[{\cite[Theorem 1.4]{ALF-ample-angles}}]
	\label{thm:polytope}
	$\overline{\h{\rm AA}(X,D)}$ is either empty or a rational polytope,
	i.e., cut out by finitely-many linear inequalities with rational coefficients in $\be_1,\ldots,\be_r$.
	\end{theorem}
	In \cite{ALF-ample-angles} we give two proofs of Theorem \ref{thm:polytope}. One of them consists on finding an explicit bijective affine transformation from a specific rational polytope to $\overline{\h{\rm AA}(X,D)}$, automatically making the latter a rational polytope. It is explicit and it requires no highbrow facts, but it is not particularly insightful. The second proof uses the fact that asymptotically log Fano varieties are Mori dream spaces. Thus, their Neron-Severi group $\mathrm{NS}(X)$ is finitely generated and the the Nef cone $\mathrm{Nef}(X)$ is polyhedral. Then, the affine map $\Phi\colon \mathbb R^r\rightarrow \mathrm{NS}(X)$ given by
	\[\Phi(\beta_1,\dots,\beta_r)= \left[-K_X-\sum_{i=1}^r(1-\beta_i)D_i\right] \]
	gives that $\overline{\hbox{\rm AA}(X,D)}=[0,1]^r\cap \Phi^{-1}({\rm Nef}(X))$ and the proof follows.
	
	The rest of \cite{ALF-ample-angles} focuses on classification. Let us recall what was known priorly. In dimension $1$, the classification is trivial: the only two possibilities are $(\mathbb P^1, p)$ and $(\mathbb P^1, p+q)$. In particular the notions of strongly ALF and ALF coincide. One notices from Example \ref{exa:ALDP.1.n} that since we may chose $n\geqslant 0$, then asymptotically log Fano varieties cannot possibly belong to a finite number of families. This is due to the nature of the definition, where $\beta$ is allowed to tend to $0$, forcing the class of these pairs to be naturally Kawamata log terminal and not $\varepsilon$-log canonical, so Birkar boundedness does not apply. Nonetheless, in dimension $2$, it is easy to prove that $X$ must be rational (\cite[\S3]{ALF}), and we can apply a version of the minimal model programme (MMP) to give a sense of classification.
	
	In \cite{ALF} the authors gave a classification of strongly ALF pairs. To do so, they introduced a notion of minimality. An ALF surface $(X,D)$ is \emph{minimal} if $X$ contains no smooth rational $(-1)$-curves $E\not\subset \mathrm{Supp}(D)$ such that $D\cdot E=1$. Moreover, they showed that either $D=D_1$ is irreducible, or $D$ consists of rational curves whose dual graph is either a chain (recall they need to intersect with simple normal crossings) or a cycle (and $D\sim -K_X$). Finally, they also showed that rational $(-1)$-curves must either be disjoint from $D$, a component of $D$, or intersect $D$ only at one component, at one point and normally.
	
	In the strong regime, \cite{ALF} noted that $(-1)$-curves within the support of $D$ can only appear at the end of chains. This simplifies the classification of strongly ALF considerably. Indeed, they apply a sort of `guided' MMP: if a surface is not minimal, then there is a $(-1)$-curve $E$ intersecting $D$ normally. They contract $E$ with a birational morphism $\pi\colon X \rightarrow \overline X$ and replace $(X,D)$ by $(\overline X, \pi_*(D)\cong D)$, which is ALF. In particular, the dual graph of $D$ and $\pi_*(D)$ are identical. One repeats this process until obtaining a minimal pair. Then, assuming $(X,D)$ is minimal, prove that either $X\cong \mathbb P^2$ or $X\cong \mathbb F_n$ for some $n\geqslant 0$. Since the intersection theory and Picard group of these curves is well understood, it is very easy to classify minimal strongly ALF pairs. In particular, any minimal strongly ALF pair $(X,D)$ has at most $4$ components in $D$. Then one can recover the rest of strongly ALF pairs by blowing up smooth points on the boundary of an ALF pair, being careful not to lose the `strong ALF' condition.
	
	\begin{remark}
		\label{rmk:weakly-ALF-difference}
	In particular, one notes that blow-ups on infinitely close points (i.e. blowing up a point on the exceptional divisor) breaks the `strong ALF' condition but not necessarily the ALF condition. 
	\end{remark}
	There is a second kind of blow-up that is `forbidden', namely blowing up two points on the (proper transform of a) fibre of $\mathbb F_n\rightarrow \mathbb P^1$. Note that the above construction (where no boundary component in $D$ is contracted) and the classification of minimal strongly ALF pairs imply that any two-dimensional strongly ALF pair (not necessarily minimal) has at most four components in its boundary \cite[Corollary 1.3]{ALF}. This makes it possible for \cite{ALF} to give a very explicit description of two-dimensional strongly ALF pairs. The case of ALF pairs is substantially more involved. Indeed, note the following fact, whose proof is left as an exercise to the reader:
	\begin{example}
		Let $(X, D=\sum_{i=1}^r D_i)$ be an ALF surface with $r\geqslant 2$. Let $p$ be a singular point of $D$ and $\pi\colon \widetilde X \rightarrow X$ be the blow-up of $p$ with exceptional divisor $E$. Denote by $\widetilde D_i$ the proper transform of $D_i$. Then $(X, E+\sum_{i=1}^r \widetilde D_i)$ is an ALF pair (but not strongly ALF). Note we replaced $D$ by its total transform.
	\end{example}
	Note the possible number of boundary components in ALF pairs is itself unbounded. Nonetheless, in \cite{ALF-ample-angles}, we were able to give a structure and a classification of minimal ALF pairs of Picard rank at most $2$.

	\begin{theorem}[{Classification of $2$-dimensional minimal ALF pairs of Picard rank at most $2$, \cite[Proposition 5.1, Corollary 5.3]{ALF-ample-angles}}]
		\label{thm:classification-rk2}
	Let $S$ be a smooth surface with $\rk(\Pic(S))\leqslant 2$, and let
	$C_1,\ldots,C_r$ be distinct irreducible smooth curves on $S$ such that
	$C=\sum_{i=1}^{r}C_i$  is a divisor with simple normal crossings.
	Then $(S,C)$ is an 
	asymptotically log Fano pair if and only if it
	is either a minimal strongly asymptotically asymptotically log Fano pair in \cite[Theorem 2.1, Theorem 3.1]{ALF} or one of the following pairs:

	\begin{itemize}
		\item Type ALdP.1.n: $S= \mathbb F_n, \, n\in\NN$, $C_1=Z_n, C_2\in |Z_n+(n+2)F|$,
		\item Type ALdP.2.n: $S= \mathbb F_n, \, n\in\NN, \, C_1=Z_n, C_2\in |Z_n+(n+1)F|, C_3\in|F|$,
		\item Type ALdP.3.n: $S= \mathbb F_n, \, n\in\NN, \,C_1=Z_n, C_2, C_3\in |F|$,
		\item Type ALdP.4.n: $S= \mathbb F_n, \, n\in\NN, \,C_1=Z_n, C_2, C_3\in|F|, C_4\in |Z_n+nF|$.
	\end{itemize}
	Moreover, the body of ample angles of the above non-strongly ALF pairs is
			\begin{equation*}
		\h{\rm AA}(S,C)=
		\begin{cases}
			\{(\be_1,\be_2)\in(0,1]^2\,:\,-n\be_1+2\be_2>0\}
			\h{ if  $(S,C)$ is $\mathrm{ALdP.1.n}$}, 
			\cr
			\{(\be_1,\be_2,\be_3)\in(0,1]^3\,:\,-n\be_1+\be_2+\be_3>0\}
			\h{ if  $(S,C)$ is $\mathrm{ALdP.2.n}$ or $\mathrm{ALdP.3.n}$}, 
			\cr
			\{(\be_1,\be_2,\be_3,\be_4)\in(0,1]^4\,:\,-n\be_1+\be_2+\be_3>0\}
			\h{ if  $(S,C)$ is $\mathrm{ALdP.4.n}$}. 
		\end{cases}
	\end{equation*}
	\end{theorem}
	Note that the acronym ALdP stands for `asymptotically log del Pezzo' which is the name for ALFs in dimension $2$.
	
		\begin{theorem}[{Structure theorem \cite[Proposition 4.1]{ALF}}]
		\label{thm:structure}
		Any two-dimensional asymptotically log Fano pair $(X,D)$ is obtained from a pair $(S,C)$ with 
		$\rk\Pic(S)\leqslant 2$
		listed in 
		Theorem \ref{thm:classification-rk2} by
		a combination of the following operations:
		\begin{enumerate}
			\item blowing-up a collection of distinct  points on 
			the smooth locus of the boundary and replacing the boundary with its proper transform,
			\item blowing-up a collection of (possibly infinitely near) singular 
			points of the boundary and replacing the boundary with its total transform.
		\end{enumerate}
	\end{theorem}
	Note that this structure result is not sufficient to provide a classification. Indeed, while its repeated application on any ALF gives at least one Picard rank $2$-model in Theorem \ref{thm:classification-rk2}, each sequence of successive blow-ups as in Theorem \ref{thm:structure} (i) from a given ALF (e.g. one in Theorem \ref{thm:classification-rk2}) will not necessarily result on an ALF. One example of obstruction is blowing-up two points in the same fibre $f\in |F|$ of the morphism $\mathbb F_n \rightarrow\mathbb P^1$, where the affected fibre is not itself in the boundary, e.g. one can do this for ALdP.4.n by blowing up two points $p_1\in C_1$ and $p_2\in C_4$ in the same component $f\in |F|$ but so that $p_1,p_2\not\in C_2\cup C_3$. The proper transform $\widetilde f$ of $f$, which is not part of the boundary, satisfies $\widetilde f^2=-2$, contradicting \cite[Lemma 2.5]{ALF}. When repeated blow-ups of type (ii) in Theorem \ref{thm:structure} take place before repeated blow-ups of type (i), or even if one intercalates them, the casuistic seems to become intractable.
	
	On the other hand, repeated blow-ups of type (ii) result in new ALFs, in so far as they are not intercalated with blow-ups of type (i). In particular, this proves that the number of components of two-dimensional ALFs is unbounded, in contrast with the strongly ALF regime. Perhaps there is a way forward to give a more explicit description of two-dimensional ALFs closer to the one for strongly ALF in \cite[Theorem 2.1, Theorem 3.1]{ALF}, but it looks unlikely due to the obstructions above.
	
	So what else can one do with ALFs beyond classifying them and study their K-stability? When it comes to classification, one can give natural generalisations of the definition to the singular setting in the sense of MMP (e.g. to canonical singularities). Whether classification is feasible is unclear, but there are questions to be answered with regards to their K-stability. Rubinstein's survey \cite{Rubinstein-survey} is a good entry point to such problems. The next natural classification problem to look into is three-dimensional ALFs. This problem remains open even when the boundary is irreducible ($r=1$), as any possible analogue of Theorem \ref{thm:structure} forces us out from the log smooth category and into the terminal one. An exception to this obstacle was noticed by Maeda \cite{Maeda} under the extra assumption that $-K_X-D$ is ample. Then, one avoids non-smooth contractions. Finally, there is a good understanding of which ALFs are K-(poly)stable, making their compactification into K-moduli components natural. However, to the best of the author's knowledge, the K-moduli of ALFs has not been explored systematically beyond some two-dimensional examples in \cite{MGPZ} where the ALF nature of the examples is really an afterthought.

		\def\bi{\bibitem}

		
		
	%
	
	
%
	{\sc University of Essex}
	
	{\tt jesus.martinez-garcia@essex.ac.uk}

\end{document}